\input amstex
\documentstyle{amsppt}
\NoBlackBoxes
\magnification 1200
\hsize 6.5truein
\vsize 8truein
\topmatter
\title Kakeya Sets and Directional Maximal Operators in the Plane \endtitle
\author  Michael Bateman  \endauthor
\affil Indiana University \endaffil
\abstract   We completely characterize the boundedness of planar
directional maximal operators on $L^p$.  More precisely, if $\Omega$ is a set of
directions, we show
that $M_{\Omega}$, the maximal operator associated to line segments in the directions
$\Omega$,
is unbounded on $L^p$, for all $p< \infty$, precisely when  $\Omega$ admits
Kakeya-type sets.  In fact, we show that if $\Omega$ does not admit Kakeya sets, then
$\Omega$
is a generalized lacunary set, and hence $M_{\Omega}$ is bounded on $L^p$, for $p>1$.

\endabstract
\thanks $email$:  mdbatema$\@ $indiana.edu 
\endthanks
\endtopmatter

\define\R{\Bbb R}
\define\om{\Omega}

\define\tom{\Cal T_{\om}}
\define\calb{\Cal B}
\define\calp{\Cal P}
\define\calt{\Cal T}
\define\cals{\Cal S}
\define\calc{\Cal C}

\define\mathcal{\Cal}

\head \S 0 Introduction \endhead

Given a closed set $\Omega \subset [0,1]$ of slopes in the plane, we let $\bold B _
{\Omega}$ be
the collection of all rectangles so that one of the sides has slope in $\Omega$, and we
define
$$ M_{\om} f (x) = \sup _{x\in R\in \bold B _ {\Omega} } {1 \over {|R|}} \int  _R f. $$
The study of such operators dates at least to Cordoba's paper [C], in which he considered
the case
$\om = [{1\over N}, {2\over N}, ..., 1]$, with the restriction that the rectangles in
$\bold
B _ {\Omega}$
have dimensions $1 \times N$.  In the case where $\om$ is a lacunary sequence, i.e., when there
is a
$\lambda \in (0,1)$ such that $\om = \{ \omega_0, \omega_1, \omega_2, ... \}$, and
$\omega_{j+1 }  \leq \lambda \omega _j$
for $j=0,1,2,...,$  Str\"{o}mberg  [S], and Cordoba and R. Fefferman [CF1] used covering
arguments
to show that $M_{\om}$ is bounded on $L^p$ when $p\geq 2$,
and Nagel, Stein, and Wainger [NSW] followed with a Fourier analytic proof for
boundedness on
$L^p$ when $p>1$.  Let us say that a set $\om$ is lacunary of order $N$ if  it is
covered by the union of a lacunary sequence $L$ of order $N-1$ with lacunary
sequences converging to every point of $L$.  Sj\"{o}gren and Sj\"{o}lin [SS] iterated the proof in
[NSW]
to improve the result to include lacunary sequences of finite order.

On the other hand, the existence of the Besicovitch set yields unboundedness of $M_{\om}$ on $L^p ,$ $p<\infty ,$ when
$\om = [0,1]$.  A further negative
result comes when $\om$ is the Cantor set:  unboundedness in this case was shown in [K]
for $p\leq 2$, and in [BK]
for $p<\infty$.

Now let us say that $\om$ admits Kakeya sets if there is a collection $\bold R_{\om}$ of
rectangles,
each pointed in a direction in $\om$ so that
$|\cup _{R \in \bold R_{\om} } R |$ is small relative to, say, $|\cup _{R \in \bold R_{\om} } 3R  |$,
where
$3R$
is the rectangle with the same center and width as $R$ and three times
the length.  In this paper we will prove
\proclaim{Theorem 0.1}  Fix $1<p<\infty$.  The following are equivalent:

A:  $M_{\om}$ is bounded on $L^p$

B:  $\om$ does not admit Kakeya sets

C:  There exist $N_1,N_2 <\infty$ such that $\om$ is covered by $N_1$ lacunary sets
of order $N_2$.
\endproclaim

To prove this theorem we will view $\om$ as being the boundary of a subtree of the binary
tree.  Then we will introduce the splitting number of a tree, which measures, loosely
speaking,  to what degree the tree has a subtree that looks like the binary tree.  This
will allow us to categorize all such $\om$ as looking like either a lacunary-type set or
a Cantor-type set.  

%We mention a few corollaries to Theorem 0.1.
%
%\proclaim{Corollary 0.2}
%The equality
%$$ \lim_{|R| \rightarrow 0}  {1 \over {|R|}} \int _R f =
%f(x) $$
%where $R$ is such that $x\in R\in \bold B _ {\Omega}$,  holds for almost every $x$
%whenever $\om$ is lacunary of finite order.
%\endproclaim
%
%
%To state the next result, we let $\Lambda \subseteq \om$ be such that $\Lambda = \{
%\omega_0 , \omega_1 , \omega_2 , ... : \omega_j \in \om \}$.  Further, we let
%$P_{\Lambda}$
%be the region of the plane bounded to the left by the line $x=1$ and bounded below by the
%line
%segments $l_k$ with slopes $-\omega _k $ connecting the points
%$ (2^k, h_k)$ and $(2^{k+1}, h_{k+1} ) $  for $k=0,1,2,...$.  Note that if we define
%$h_0= 1$,
%then the slopes $\omega _k$
%will determine
%the other values $h_k$.  With this notation, the Theorem of A. Cordoba and R. Fefferman
%given in [CF] implies this result:
%
%\proclaim{Corollary 0.3}
%If we define
%$$ \widehat{T_{\Lambda} f }= \chi _{P_{\Lambda} } \hat{f}, $$
%then $T_{\Lambda}$ is bounded on $L^p, 1<p<\infty$, if $\om$ is lacunary of finite order.
%\endproclaim

\proclaim{Acknowledgements}
The author thanks Russ Lyons, and especially Nets Katz, for helpful discussions.
\endproclaim

\head \S ${1\over 2}$ Outline \endhead

The goal of this paper is to provide a proof of Theorem 0.1.  The proof that
$A\Rightarrow B$ is simple.
For suppose $\Omega$ admits Kakeya sets in the sense above:  then for any $N$, we have
sets
$$E_N : = \bigcup _ t R_t ^{(N)}  \text{ \quad and\quad } E_N ^* : = \bigcup _ t 3R_t
^{(N)}, $$
where the slopes of the $R_t ^{(N)}$ are in $\Omega$, such that
$$  {    {| E_N | } \over {| E_N ^*| }   }   \rightarrow 0 \text {\quad  as \quad }
N\rightarrow \infty ,$$
and such that $M_{\om} \chi _{E_N} (x) > {1\over 2} $ when $x\in E_N ^*$.  Hence
$${    {     \int _ {\R ^n }  ( M_{\om}  \chi _{E_N}  ) ^p   } \over
{\int _ {\R ^n }  ( \chi _{E_N} ) ^p   }    }   \gtrsim {    {| E_N | } \over {| E_N ^*| }
 } \rightarrow \infty \text {\quad  as \quad }
N\rightarrow \infty ,  $$
where, of course, by $\alpha \lesssim \beta$ we mean $\alpha \leq c \beta$ for some
constant $c$.

Our contribution is the proof that $B\Rightarrow C$, and the majority of the paper is
devoted to this.  For completeness, we will review in the final section a proof that $C\Rightarrow A$.

\head \S 1 Splitting Number and Notation for Trees \endhead

%We begin by reviewing some notation for trees and then introducing the splitting number.
%All trees considered in this paper will be subtrees of the binary tree.  More precisely:
%We define the \it height \rm of a vertex $v$, $h(v)$, to be the length of the shortest path to the
%origin.  (Here we will use ``origin" instead of the more commonly used ``root".)  We say
%that
%$u$ is a \it child \rm of $v$ if there is an edge connecting $u$ and $v$ and $h(v)=h(u)-1$.  So by
%the binary tree $\Cal B$, we mean the tree such that every vertex has exactly two
%children.  

We begin by constructing the binary tree $ \calb  $:  fix a vertex, $v_{origin}$, called the 
\it origin, \rm and define 
$$\calb _ 0 = \{ v_{origin} \}.$$
(Here we will use the word ``origin", since the
more commonly used word ``root" will be used frequently as a verb.)  Then for $n=0,1,2,...$, 
connect each vertex $v\in \calb _n$ to two new vertices $c_0 (v)$ and $c_1 (v)$, called the $0$th and $1$st
children of $v$,  and define
$$\calb _{n+1} = \bigcup _{v\in \calb _n } \{c_0 (v) , c_1 (v) \} .$$
Then $\tilde{ \calb } $ is the tree with vertices 
$$\calb : = \bigcup _{n=0} ^{\infty} \calb _n $$
and edges connecting a vertex $v$ with its children $c_0 (v)$ and $ c_1 (v) $.  We will refer to the tree 
$\tilde{ \calb } $ by its vertex set $\calb$, and we will do the same for other trees considered in this
paper, which will all be subtrees of $ \calb  $.  
If $v\in \calb _n$, define the height of $v$, $h(v) = n$.  Further, if $\calt \subseteq \calb $, then by $\calt _ k$ we mean all vertices $v\in
\calt $ such that
$h(v) = k$.  

Now given a vertex $v\in \mathcal T \subseteq
\calb$, we define a \it ray $R$ rooted at $v$ \rm to be an ordered set of
vertices $(v_1 = v, v_2, v_3, ...)$ such that $v_{j+1}$ is a child of $v_j$ for $j =
1,2,...$  Loosely speaking, a ray rooted at $v$ is a path from $v$ to infinity that
always moves (strictly) away from the origin of the tree.  The \it boundary  \rm  of $\calt$ is the set
of all rays in $\calt$ rooted at the origin, and will be denoted $\partial \calt$.  
Define the \it shadow, \rm  $U(v)$, of a vertex $v$,  to be the set of all rays $R$ such that $v\in R$.  

%Similarly, if $\om$ is any closed subset of $[0,1]$, we will let $\mathcal
%T_{\om}$ be the subset of $\calb$ formed by the rays in the boundary identified with the
%numbers  $.a_1a_2...$ in $\om$.
%(Note that this can be done since if $\om$ is closed, then $[0,1] - \om$ is a collection
%of open intervals, each of which can be covered by a union of dyadic intervals.  Then
%$\mathcal T_{\om}$ will be $\mathcal T_{[0,1]} = \calb$ minus all rays that are identified
%with numbers in those dyadic intervals.)

We identify 
the vertices of the binary tree with the dyadic intervals contained in $[0,1]$ as follows:
$$\eqalign{
&1.  \text{  Identify the origin with $[0,1]$.  } \cr
&2.  \text{  If $v$ is identified with the dyadic interval $I$,   } \cr
&  \text{ then identify $c_0 (v) $ with the left half of $I$, } \cr
&  \text{  and identify $c_1 (v) $ with the right half of $I$.   }    }$$
We will write $v_I$ to indicated the vertex identified with the interval $I$, and $I_v$ to indicate the interval 
identified with the vertex $v$.  
Similarly, we can identify the boundary of the binary tree with the interval
$[0,1]$ in the following natural way:  identify $.a_1a_2...$, where $a_j \in \{ 0,1 \} $,
with the ray $(v_0 = v_{origin}, v_1, v_2, ...) $ if $v_{j+1}  $ is the $a_{j+1}$th child of $v_j$,
i.e., if $v_{j+1} = c_{a_{j+1}} (v_j) $ for every $j= 1,2,...$.

If $\Omega$ is closed, then $[0,1] - \Omega$ is the union of open intervals $Q_j$.  Each $Q_j$ is 
the union of dyadic intervals, so we may write 
$$ [0,1] - \Omega = \bigcup _j I_j ,$$
where each $I_j$ is a dyadic interval.  We define $ \calt _{\Omega} $ to be the subtree of $\calb$ 
obtained by removing the subtrees of $\calb$ rooted at $v_{I_j} $ for $j=1,2,...$.  Alternatively, $\calt _{\Omega} $ 
is the subtree of $\calb$ with boundary
$$\partial \calt _{\Omega} = \partial \calb - \left( \bigcup _{j=1} ^{\infty} U (v_{I_j} ) \right) . $$

Earlier, we defined what it means for a ray $R$ to be rooted at a vertex $v$. 
The collection of such rays depends on the tree $\calt$, and will be denoted by $\goth R _{\mathcal T}  (v)$.
If $u \in R$ for some $\goth R _{\mathcal T}  (v)$, we will write $u \subseteq v$, and
say that
$u$ is a \it descendant \rm  of $v$, or that $v$ is an \it ancestor \rm of $u$.

We will say that a vertex $v$ \it splits \rm, or call $v$ a \it splitting
vertex \rm, if
$v$ has two children, and define the splitting number split$(R)$ of a ray $R$ to be the number
of splitting vertices along $R$.  Then the splitting number of a vertex $v$ with respect
to a tree $\cals$ rooted at $v$ is defined to be 
$$split_{\mathcal S} (v) = \min _{R\in  \goth R _{\cals}  (v) } split (R)$$
and the splitting number of $v$ is defined to be 
$$split (v) = \sup _ { \mathcal S  }  split_{\mathcal S} (v),$$
where the sup is taken over all subtrees $\mathcal S$ of $\mathcal T$ rooted at $v$.
Finally, for a tree $\mathcal T$, we define
$$split(\mathcal T) = \sup _{v\in \mathcal T } split (v).$$
Before we state a theorem using this new language, we give a definition of lacunarity
that is more suitable for trees:  a subtree $\mathcal L \subseteq \mathcal T$ is said to
be
lacunary of order $0$ if $\mathcal L$ consists of a single ray in the boundary of
$\mathcal T$, and $\mathcal L$ is said to be lacunary of order $N$ if all splitting
vertices of $\mathcal L$ lie along a lacunary tree of order $N-1$.
%\proclaim{Remark 1.1}
%We will show in \S $5$ that the boundary of a lacunary tree of order $1$ can be covered by four lacunary sequences in
%the traditional sense.
%\endproclaim

\proclaim{Remark 1.1}
If $\mathcal L$ is a lacunary tree of order $1$, then, 
loosely speaking, the directions associated with $ \mathcal L $ can be covered by four lacunary 
sequences in the traditional sense.  More precisely, define $\alpha (R)$ to be the real number 
in $[0,1]$ identified with the ray $R$.  Note that all splitting vertices of $\mathcal L$ lie on a single ray, call it $ \lim \mathcal L  $. 
We claim there exist four lacunary sequences
 $\{a_j ^{(i)} \} _{j=1} ^{\infty}$,  $ i =1,..., 4$, such that  
$$a_j ^{(i)}  \rightarrow \alpha (\lim \mathcal L ) \quad \text{as} \quad j \rightarrow \infty $$ 
for each $i$, and such that 
$$\alpha (R) \in \left( \bigcup _{i,j} a_j ^{(i)}  \right) \cup   \{  \alpha (\lim \mathcal L ) \}  $$
for every $R \in \partial \mathcal L$.  
\endproclaim
\demo{Proof}
For each $j=0,1,2,...$ there is 
at most one $R \in  \partial \mathcal L$ such that $d(\alpha (R),  \alpha (\lim \mathcal L ) ) = 2^{-j}$, where $d$ 
denotes the dyadic distance on real numbers.  (That is, $d(\beta_1 , \beta_2 )$ is defined to be the size of the smallest dyadic 
interval containing both $\beta_1$ and $\beta _2$.)  Now consider the set
$$ A: = \{ R \in \partial \mathcal L \colon \alpha (R) >  \alpha (\lim \mathcal L )  \} .$$
Finally, observe that if $R_0, R_1, R_2, ... \in A$ are such that 
$$d(\alpha (R_j),  \alpha (\lim \mathcal L ) ) = 2^{-2j} , $$ then 
$$ 0 < \alpha ( R_{j+1} ) -  \alpha (\lim \mathcal L ) \leq {1 \over 2}  
\left( \alpha ( R_{j} ) -  \alpha (\lim \mathcal L ) \right) ,$$
i.e., $\{\alpha (R_j) \} _{j=1} ^{\infty} $ is lacunary in the traditional sense.  An identical claim can be made if 
$R_0, R_1, R_2, ... \in A$ are such that $d(\alpha (R_j),  \alpha (\lim \mathcal L ) ) = 2^{-2j-1}$, hence
$A$ is covered by two lacunary sequences in the sense described above.  Of course this implies $\mathcal L$ 
is covered by four lacunary sequences since we could similarly show that the set 
$B : = \{ R \in \partial \mathcal L \colon \alpha (R) <  \alpha (\lim \mathcal L )  \} $ is covered by two 
lacunary sequences. \qed
\enddemo

\proclaim{Theorem 1.2}

A:  If split$(\mathcal T_{\om}) = N < \infty$, then $\mathcal T_{\om} $ is lacunary of
order $N$, and
hence $M_{\om}$ is bounded on $L^p$ for $1<p<\infty$.

B:  Conversely, if split$(\mathcal T_{\om}) =  \infty$, then $\om$ admits Kakeya sets,
and hence $M_{\om}$ is unbounded on $L^p, p<\infty$.
\endproclaim
\proclaim{Remark 1.3}
Let $\Omega$ be such that split$(\mathcal T_{\om}) = N$.  In \S $5$ we will see that there 
exists a constant $C$ such that 
$$ || M_{\Omega} f|| _p \leq C N || f ||_p .$$

%Let $\Omega$ such that split$(\mathcal T_{\om}) = N$. In light of Remark 1.1 and Theorem
%1.2,
%we have that if $L = \{  { 1\over 2 } , { 1\over 4 }, { 1\over 8 }, ... \} $,
%and if $ || M_L f|| _p \leq A_p ||f||p$ for $f\in L^ p$,
%then
%$$ || M_{\om} f|| _p \leq 4^N  N A_p ||f||p \text{\quad for  \quad } f\in L^ p .$$
\endproclaim
Theorem 1.2 automatically yields the $``B\Rightarrow C"$ part of Theorem 0.1, and
we are already able to dispense with part A of Theorem 1.2.  The following lemma records
an easy observation that will help with the proof.
\proclaim{Lemma 1.4}
If $\mathcal T$ is a tree, and $u\neq v$ are vertices of $\mathcal T$ with split$(u) \geq
N$, and split$(v) \geq N$, and if $h(u) \geq h(v)$, then either split($\mathcal T) \geq N+1$, or there exists
$R\in \goth R _{\mathcal T}  (v)$ such that $u \in R$, i.e., $u\subseteq v$.
\endproclaim
\demo{Proof}
First note that $u$ and $v$ must have a common ancestor.  If there is no $R\in \mathcal R
_{\mathcal T}  (v)$ such that $u \in R$, then the common ancestor is some other vertex
$w$, and
$v\neq w\neq u$.  (Of course $u$ cannot be the common ancestor since $h(u) \geq h(v)$.)
But then split$(w) \geq N+1$: since there are subtrees $\mathcal T_v$ and  $\mathcal T_u$
for which split$_ {\mathcal T_v }  (v ) = N =$ split$_ {\mathcal T_u }  (u )$, we define
$\mathcal T_w$
to be the tree formed by joining $\mathcal T_u$ with $\mathcal T_v$ through $w$, 
and we have split$_{\calt _w} (w) = N+1$. \qed
\enddemo
\demo{Proof of Theorem 1.2 part A}
If split$(\mathcal T_{\om})    = 0$, then $\mathcal T_{\om}$ has only one ray rooted at
the origin.  Hence
$\mathcal T_{\om}$ is
lacunary of order zero.  Now we induct on the splitting number:  suppose split$(\mathcal
T_{\om}) =N$.
By Lemma 1.4, all vertices $v$ with split$(v)=N$ lie along a single ray.  So if $v^*$ is
a child of $v$ not lying on the ray $R$, and if $\mathcal T _{v^*}$ is a subtree of
$\mathcal
T_{\om} $ rooted at $v^*$, then split$_ {\mathcal T_{v^*} }  (v^* )\leq N-1$,
and hence is lacunary of order $N-1$ by the induction hypothesis.
But we can repeat this process, which results in $T_{\om}$ being lacunary of order $N$. \qed
%But then $\mathcal T_{\om}$ consists only of lacunary trees of order $N-1$ along the rays
%of a lacunary sequence of order $1$, so $\mathcal T_{\om}$ is lacunary of order $N$. \qed
\enddemo

Since we suppose now that split$(\mathcal T_{\om} ) = \infty$, to prove Theorem 1.2 part
B it suffices
to exhibit, when split$(\mathcal T_{\om} )  \geq N$, a collection of parallelograms $\{
P_t
\}$, each of which is pointed in one of the directions in $\om$, such that
$$ |\cup _t P_t | \lesssim {1\over N}     \tag{$ \clubsuit $ 1} $$
and such that
$$ |\cup _t 3P_t | \gtrsim {{\log N} \over N}, \tag {$ \clubsuit $ 2} $$
where $3P_t $ is the parallelogram with the same center and width as $P_t$ , but three
times the length.  To do this, we divide the interval $[0,1]$ on the $y$-axis into small
intervals, each of which wil serve as a base for one of the parallelograms $P_t$.  The
difficult part of the construction is to specify a slope for each of the parallelograms
so that they satisfy the properties
$\clubsuit 1$ and $\clubsuit 2$.  In fact, we will not give an explicit choice of slopes;
instead, we will use the probabilistic method to show that such a choice exists.

\head \S 2 Pruned Trees and Sticky Maps \endhead

It will actually be to our advantage to limit the possible slopes to a subset of $\om$,
represented by a $pruned$ subtree $\calp$ of $\calt _{\om}$, and to restrict our
attention to
a certain class of slope
functions, called \it sticky \rm maps.  We now define these terms.  Suppose $\mathcal T
\subseteq
\mathcal B$ is a tree such that split$(\mathcal T) =N$.  Then there is a vertex $v_0 \in
\mathcal T$ such that split$(v_0) =N$.  Without loss of generality, suppose $v_0$ is the
origin.  We say that $\mathcal T$ is $pruned$ if for every $R\in \goth R _{\mathcal T}
(v_0)$, and every $j=1,2, ..., N$, $R$ contains exactly one vertex $v_j$ such that
split$(v_j) =j$.  If $\mathcal T$ is not necessarily pruned, then we can find a pruned
subtree $\mathcal P$ of $\mathcal T$ that still has splitting number $N$ by the following
recursive procedure:

$$\eqalign{
&1.  \text{  Let $v_0$ (the origin) be in $\mathcal P$.  } \cr
&2.  \text{  Assign $j: = 0$.   } \cr
&3.   \text{  While $0 \leq j < N$, if $v\in \mathcal P$ has splitting number $N-j$, }
\cr
&\text{ \quad      choose a pair $u,w$ such that $u\subseteq c_0 (v)$,  $w\subseteq c_1
(v)  $,  } \cr
% &\text{ \quad        } \cr
&\text{ \quad       split$(u) \geq N-j -1$, and split$(w) \geq N-j -1,$ and add $u,w$ to
$\calp$. } \cr
&\text{ \quad      Also add all vertices and edges connecting $v$ to $u$ and $w$. } \cr
%&\text{ \quad       if $u,w \in \mathcal T$ have splitting number $j-1$, } \cr
%&\text{  \quad       if $u,w$ are descendants of $v$ in $\mathcal T$, } \cr
%&\text{  \quad       if $u$ is a descendant of the $0$th child of $v$   } \cr
%&\text{  \quad       and if $w$ is a descendant of the$1$st child of $v$, } \cr
%&\text{  \quad       then add $u$ and $w$ to $\mathcal P$.  }  \cr
%(And, of course, add the edges connecting $u$ and $w$ to $v$. )   } \cr
&4.  \text{  Assign $j := j+1$.  }
}$$

We call the vertices added to $\calp $ at the $j$th iteration the $j$th generation, and
denote the
collection of vertices in the $j$th generation $G_j ( \mathcal T)$.  If $\calt$ is
already pruned, then
$G_j (\calt )$ still makes sense.  We will denote by $\calp (\calt )$ the subtree of
$\calt$ formed by
$\cup _j G_j (\calt ) $ and the edges and vertices connecting $G_j (\calt ) $ to $G_{j+1}
(\calt ) $.  Note that it is not
necessary for any $v\in G_j$ to have $h(v)=j$.  (Except for $j=0$, because we have
assumed $h(v_0)=0$, and $G_0 = \{ v_0 \}.$)  Also note that for a general tree $\mathcal
T$ with splitting number $N$, there may exist several different subtrees $\mathcal P _1,
\mathcal P _2, \mathcal P _3, ...$ each with splitting number $N$, and each pruned.  The
method above yields one of them.

We now consider maps $\sigma : \calb \rightarrow \mathcal S \subseteq \calb$.  Such a map
is said to be $sticky$ if whenever $u \subseteq v \in \calb$, then $\sigma (u)  \subseteq
\sigma (v) $ in $ \mathcal S$.  In addition, all sticky maps considered here will be assumed
to satisfy $h(\sigma (v)) = h(v)$ for all $v\in \calb$.   Recall that $I_v$ is the dyadic interval 
identified with $v$, and note that $| I_v | = 2^{-h(v)}$, where $| \cdot |$ denotes the standard 
Euclidean measure.  Note that if $V$ is a collection of vertices, and if
$v_1, v_2, ...$ are the disjoint maximal elements in $V$, then 
$$ |\cup _{v\in V} I_v | =\sum | I_{v_j} | . $$  
The following lemma gives an elementary fact about sticky maps into pruned trees
that will be useful in completing the proof of Theorem 1.2 part B.
\proclaim{Lemma 2.1}
Suppose $\sigma : \mathcal B \rightarrow \mathcal P$ is a sticky map into a pruned tree $\calp$
with generations $G_k (\calp ) $.   
Then for each $k = 0,1,..., N$, 
$$ \sum _{v\in G_k (\calp ) } \sum _{u \in \sigma ^{-1} (v) } |I_u|   =1. $$
\endproclaim
\demo{Proof}
Of course, for $v_1 \neq v_2 \in  G_k (\calp )$, the sets $\sigma ^{-1} (v_1)$ and $\sigma ^{-1} (v_2)$ 
are disjoint.  Hence the sum in the statement of the lemma is over a collection of vertices 
$u_1, u_2, ... \in \calb$ such that $u_i \not\subseteq u_l$ whenever $i\neq l$.  But then 
$$ \sum _{v\in G_k (\calp ) } \sum _{u \in \sigma ^{-1} (v) } |I_u|    = 
 \left| \bigcup _{ u \in  \sigma ^{-1} ( G_k (\calp )  )  } I_u \right| \leq 1 ,$$
since $|[0,1] | = 1$.  In fact, equality holds, since 
%if $t\in \calb _{ h(\calp ) }$, then there is some 
%$u \in  \sigma ^{-1} ( G_k (\calp )  ) $ such that $t\subseteq u$.  Hence 
$ \bigcup _{ u \in  \sigma ^{-1} ( G_k (\calp )  )  } I_u = [0,1]$.
\enddemo

%\proclaim{Lemma 2.1}
%Suppose $\sigma : \mathcal B \rightarrow \mathcal P$ is a sticky map into a pruned tree
%$\calp$ with splitting number $N$.
%
%A:  If $v_2 \subsetneqq v_1$, and if $v_1 \in \sigma ^{-1} (G_k (\calp) )$, then $v_2
%\notin \sigma ^{-1} (G_k (\calp) )$.
%
%B:  Let $|\cdot|$ denote the measure on $\calb$ defined above, and let $G_k (\calp)$ be
%the $k$th generation of $\calp$.  Then
%$$ \sum _{u \in G_k (\calp) }    \# (\sigma ^{-1} (u) )   |u|  = 1.$$
%\endproclaim
%\demo{Proof}
%If $v_1, v_2 \in G_k (\calp )$, then $v_2 \not\subseteq v_1$.
%But $\sigma ( v_2 ) \subseteq \sigma ( v_1 )$ since $\sigma$ is sticky.  This proves part
%A.
%\enddemo

\head \S 3 Geometric Construction \endhead

We can now be more specific about how to construct the collection of parallelograms
mentioned above.  We define the height of a tree $\calt$ by
$$h(\calt) = 1+ \sup _{v\text{  splitting} }  h(v)$$
where the sup is only taken over vertices $v$ that split.  Since we suppose
split$(\mathcal T_{\om})
\geq N,$ there exists a pruned subtree $\calp : = \calp ( \calt _{\om}  )  $ such that
split$(\calp) = N$.
We will ignore all vertices $v\in \calp$ such that $h(v) > h(\calp )$.  
For each $t = 0, {1 \over {  2^{h(\calp) }  }  }, {2 \over {  2^{h(\calp) }  }  }, ... ,{
{  2^{h(\calp) } -1 } \over {  2^{h(\calp) }  }  },  $ we will have a parallelogram
$P_t =P_{t,\sigma}$
with corners $(0,t), (0, t+ 2^{h(\calp) } )  , (2, t+ 2\sigma (t) ) $, and $(2, t+
2^{h(\calp) } + 2\sigma (t) )  $, where $\sigma : \calb \rightarrow \calp$ is a sticky map
to be determined.  We will write
$$
K_{\sigma} = \bigcup _t P_{t,\sigma}.
$$

To finish the proof of Theorem 1.2 part B, it remains to prove the
following.
\proclaim{Claim 3.1}
A:  If $\sigma : \calb \rightarrow \calp$ is sticky, then
$$ | K_{\sigma} \cap \left( [0,1] \times \Bbb R \right) |
\gtrsim {   {\log N } \over {N } }.$$

B:  There exists a sticky map $\sigma : \calb \rightarrow \calp$ such that
$$  |K_{\sigma} \cap \left( [1,2] \times \Bbb R \right) |
\lesssim {1\over N}.$$
\endproclaim
%\demo{Proof of A}
We begin by proving A.  For $j=0, 1,2,...$, define $X_j$ to be the vertical strip
$$X_j = [2^{-j}, 2^{-j+1}] \times \R.$$
We will show that for $j=0,1,..., \log N$, we have the estimate
$$ | K_{\sigma} \cap X_j | \gtrsim {   1 \over N  },$$
and Claim 3.1 A will follow.  To do this we will control the intersections of the
rectangles $P_{t, \sigma }$ in the strip $X_j$.  A more precise statement is given below
in the
setting of a measure space:
\proclaim{Lemma 3.2}
Let $(X, \mathcal M , |\cdot | )$ be a measure space, and let $A_1, A_2, ..., A_K$ be sets
with $|A_i | = \alpha$ for every $i$.  If
$$ \sum _{i = 1} ^ K   \sum _{l = 1} ^ K|  A_i \cap A_l | \leq M, $$
then
$$ |\bigcup _{i=1} ^K  A_i | \geq {  {\alpha ^2 K^2 } \over { 16M} }.$$
\endproclaim
\demo{Proof of Lemma 3.2}
By pigeonholing, we obtain a set $E\subseteq {1,..., K}$ such that $\# (E)  \geq {K\over
2}$, and
$$ \sum _{l = 1} ^ K | A_i \cap A_l | \leq  {{2M} \over K} $$
for $i\in E$.  But this implies that
$$ {1\over {\alpha } } \int _{A _i} \sum _{l=1 } ^ K \chi _{A_l} \leq { {2M} \over
{\alpha K} }$$\
for $i\in E$, and hence that
$$ \sum _{l=1 } ^ K \chi _{A_l} (x) \leq { {4M}\over { \alpha K } }$$
for $x$ in a set $B_i \subseteq A_i$, with $B_i \geq { {\alpha} \over 2}$ for $i\in E$. 
But then
$$ | \bigcup _{i=1} ^K  A_i |  \geq  |\bigcup _{i\in E} ^K  B_i | \geq      { {\alpha K }
\over {4M } }    \sum _{i\in E } | B_i| = {  {\alpha ^2 K^2 } \over { 16M} }. $$
\qed
% End of proof of Lemma 3.2
\enddemo
Writing $P_{t,\sigma , j } : =  P_{t,\sigma } \cap X_j$, we have
$$| P_{t,\sigma , j } |= 2^{-h(\calp ) - j }. $$
Then by Lemma 3.2, we only need to show
$$ \sum _{t_1 = 1} ^ {2^{h(\calp ) } }  \sum _{t_2 = 1} ^ {2^{h(\calp ) }  }   |
P_{t_1,\sigma , j } \cap P_{t_2,\sigma , j } | \lesssim {  {N} \over {2^{2j} } }.
\tag{$ \pitchfork $}   $$
Since the diagonal term is
$$ \sum _{t_1 = 1} ^ {2^{h(\calp ) } }   |  P_{t_1,\sigma , j } | = {1\over {2^j } }, $$
we will only be able to show ($ \pitchfork $) when $j= 0,1, ..., \log N$.  Let us introduce some notation that will
be helpful in decomposing the sum in ($ \pitchfork $).  For any two vertices $u$ and $v$ in $\calb _{h ( \calp ) }$, 
let $D (u,v)$ be the minimal vertex containing both $u$ and $v$, i.e., let $D (u,v)$ be the 
vertex $w$ with largest height
satisfying $u\subseteq w$ and $v\subseteq w$.  (Of course this notion could be defined on all pairs of vertices
in $\calb$, but we only need to use the restriction to pairs of vertices in $\calb _{h ( \calp ) }$.)  Then of course 
for a vertex $w$, we will write
$$ D ^{-1} (w) : = \{ (u,v) \in \calb _{ h(\calp ) }   \times \calb _{ h(\calp ) } \quad \colon \quad D (u,v) =w \},  $$
and if $W$ is a collection of vertices, we will write $D^{-1} (W) : = \cup _{w\in W} D^{-1} (w) $.    Now suppose $t_1 \neq t_2$, and note that if
$$ P_{t_1,\sigma , j } \cap P_{t_2,\sigma , j } \neq \emptyset ,$$
then
$$ 2^{-j} |I_ {D(t_1, t_2 ) } | \gtrsim 2^{-j} | \sigma (t_1) - \sigma (t_2) | \gtrsim |t_1 - t_2|.  \tag{$ \diamondsuit $ } $$
(Note that $|I_ {D(t_1, t_2 ) } |$ is just the usual dyadic distance, except that here it is defined on vertices in $\calb$ 
identified with $ {h ( \calp ) } -$digit binary expansions.)
In light of this, we introduce for a vertex $w\in \calb$,
$$ \Gamma ^j (w) : = \{ (t_1 , t_2 ) \in D ^{-1} (w) \quad \colon \quad 2^j |t_1 - t_2 | \lesssim |I_w| \} $$
and
$$ \Gamma ^j _l (w) : = \{ (t_1 , t_2 ) \in D ^{-1} (w) \quad \colon \quad 2^{j+l} |t_1 - t_2 | \sim |I_w| \} $$
so that
$$ \Gamma ^j (w) = \bigcup _{l\geq 0}  \Gamma ^j _l (w) . $$
Observe that 
$$ \#  \left( \Gamma ^j _l (w) \right)  \sim 2^{2h( \calp ) -2l -2 j - 2h(w) }  $$
and 
$$ \#  \left( \Gamma ^j  (w) \right)  \sim 2^{2h( \calp )  -2 j - 2h(w) } . $$

Now we write the off-diagonal part of ($ \pitchfork $) as 
$$\eqalign{
\sum _{t_1 \in \calb _{ h (\calp ) }    }  \sum _{t_2 \neq t_1}
| P_{t_1,\sigma , j } \cap P_{t_2,\sigma , j } | 
& = \sum _{v\in \calp } \sum _{ (t_1 , t_2 ) \in D^{-1} (\sigma ^{-1} (v) ) }
| P_{t_1,\sigma , j } \cap P_{t_2,\sigma , j } |   \cr 
& = \sum _{v\in \calp } \sum _{u\in \sigma ^{-1} (v) }
\sum _ { (t_1 , t_2 ) \in D^{-1} (u) }  
| P_{t_1,\sigma , j } \cap P_{t_2,\sigma , j } |  \cr
&=:   ( \star ) .  } $$

%& \lesssim \sum _{k=1} ^N \sum _{v\in G_k (\calp )  }  \sum _{u\in \sigma ^{-1} (v) }
%\sum _ { (t_1 , t_2 ) \in D^{-1} (u) }  
%| P_{t_1,\sigma , j } \cap P_{t_2,\sigma , j } |  \cr
%&=:   (\star)   } $$

By ($\diamondsuit$), we have that if $(t_1 , t_2 ) \in D^{-1} (u)$ is to contribute to the sum, then 
$$ |I_u| \gtrsim 2^j |t_1 - t_2 |,   \tag{  $\diamondsuit \diamondsuit $ }  $$
i.e., $(t_1, t_2 ) \in  \Gamma ^j   (u) $.  
Further, if $(t_1 , t_2 ) \in \Gamma ^j _l  (u)$ contributes, then 
$$ |\sigma ( t_1 ) - \sigma (t_2 ) | \gtrsim 2^j | t_1 - t_2 | ,$$
so that in this case,
$$ | P_{t_1,\sigma , j } \cap P_{t_2,\sigma , j } |  \lesssim { 1   \over { 2^{ 2 h(\calp ) } 2^j |t_1 - t_2 |  } } 
\sim    { 1   \over { 2^{ 2 h(\calp ) } 2^{-l} |I_u |  } }  . $$
Then because of ($\diamondsuit \diamondsuit $), we may compute the innermost sum in ($\star $) as
$$ \eqalign{   \sum _ { (t_1 , t_2 ) \in D^{-1} (u) }  
| P_{t_1,\sigma , j } \cap P_{t_2,\sigma , j } |
& =    \sum _ { (t_1 , t_2 ) \in \Gamma ^j  (u) }  
| P_{t_1,\sigma , j } \cap P_{t_2,\sigma , j } |  \cr 
& = \sum _{l\geq 0}  \sum _ { (t_1 , t_2 ) \in \Gamma ^j _l  (u) }  
| P_{t_1,\sigma , j } \cap P_{t_2,\sigma , j } |  \cr 
& \lesssim  \sum _{l\geq 0}  {    {  \# \left(   \Gamma ^j _l  (u)    \right)   } \over {  2^{2h(\calp)  - l } |I_u |  } }   \cr
& \lesssim { { |I_u| }  \over  { 2^{2j}  }   }  .    }$$

To finish estimating $(\star)$, we state and prove a technical-looking proposition, whose proof requires little 
more than counting exponents.

\proclaim{Proposition 3.4}
Fix $w\in \calb $.  Then for any $l^* \geq 1 $,
$$ \sum _{l=0} ^{l^* -1 } \sum _{\{ u\subseteq w :h(u) = h(w) +l \} }  \# \left( \Gamma ^j _{l^* - l} (u) \right)  
\leq 2   \sum _{ \{ u\subseteq w : h(u) = h(w) +l^*  \} }  \# \left( \Gamma ^j (u) \right) . $$
\endproclaim
\demo{Proof}
There are $2^l$ vertices $u\subseteq w$ such that $h(u) = h(w) + l$, so the estimates on $\Gamma ^j _l (w)$ 
and $\Gamma ^j  (w)$ allow us to control the left hand side in the statement of the proposition by 
$$ \sum _{l=0} ^{l^* -1 } 2^l 2 ^{2h( \calp ) -2(l^* -l ) -2 j - 2(h(w) + l ) } \lesssim 2^{l^*} 2^{2h( \calp ) -2 j - 2h(w) }, $$
which is controlled by the right hand side. \qed
\enddemo
This Proposition allows us to restrict attention in the outer sum in $(\star)$ to splitting vertices $v \in \calp$, i.e.,
to vertices $v \in G_k (\calp ) $ for some $k = 1,2,..., N$.  
For if $v_1, v_2, ... \in \calp$ are such that $v_{j+1}$ is a child of $v_j$, with 
$v_1, v_2, ..., v_{n-1} $ not splitting and $v_n$ splitting, then by Proposition 3.4, 
$$\sum _{l=1} ^{n-1}   \sum _{u\in \sigma ^{-1} (v_l ) } 
\sum _ { (t_1 , t_2 ) \in D^{-1} (u) }  | P_{t_1,\sigma , j } \cap P_{t_2,\sigma , j } |
\lesssim   \sum _{u\in \sigma ^{-1} (v_n ) }   { { |I_u| }  \over  { 2^{2j}  }   } .$$
%\sum _ { (t_1 , t_2 ) \in D^{-1} (u) }  | P_{t_1,\sigma , j } \cap P_{t_2,\sigma , j } |  . $$
Using the computation above with Proposition 3.4 and Lemma 2.1 yields
$$\eqalign{   (\star )  & \lesssim  
%\sum _{k=1} ^N \sum _{v\in G_k (\calp )  }  \sum _{u\in \sigma ^{-1} (v) }
%\sum _ { (t_1 , t_2 ) \in D^{-1} (u) }  
%| P_{t_1,\sigma , j } \cap P_{t_2,\sigma , j } |  \cr
 \sum _{k=1} ^N \sum _{v\in G_k (\calp )  }  \sum _{u\in \sigma ^{-1} (v)  }   
{ { |I_u| }  \over  { 2^{2j}  }   }  \cr 
& \lesssim {N \over { 2^{2j}  }   },
}$$
which completes the proof of Claim 3.1 A.

\head \S 4 The Probabilistic Argument and Percolation on Trees \endhead

Now we prove probabilistically that there is some sticky map $\sigma  : \calb \rightarrow
\calp$ such that $ | K_{\sigma } \cap \left( [1,2] \times \Bbb
R \right) | \lesssim {1\over N}$.  In fact, if we denote by $Pr(x,y)$ the probability
(over sticky $\sigma$) that $(x,y) \in P_{t,\sigma}$ for some $t$, it is enough to show
that given
$(x,y) \in[1,2]\times[0,3]$,
we have $Pr(x,y) \lesssim {1\over N}$.  (Of course, if $x\in [1,2]$, and $y\notin [0,3]$,
then $(x,y)$ cannot possibly be covered by $K_{ \sigma}$.)  Then by the linearity of
expectations, we would have
$$\eqalign{     \int \left( \int _1 ^2 \int _0 ^3 \chi _{ K_{\sigma } } (x,y)
dydx \right) d\sigma  & =   \int _1 ^2 \int _0 ^3  \left( \int \chi _{ K_{\sigma } }
(x,y) d\sigma \right) dydx \cr
& =  \int _1 ^2 \int _0 ^3 Pr(x,y) dydx  \cr
& \lesssim {1\over N}.  }$$
This, of course, would imply the existence of a sticky map $\sigma$ satisfying Claim 3.1
B.
%So now we construct sticky maps $\sigma  : \calb \rightarrow \calp$ randomly as
%follows:

%$$\eqalign{& 1.  \text{  sSend the origin of $\calb$ to the origin of $\calp$. } \cr
%&2.   \text{  If $\sigma  : \calb \rightarrow \calp$ does not split, }  \cr
%&\text{ \quad and if $u$ is a child of $v$,  } \cr
%&\text{   \quad     then set $\sigma (u)$ equal to the child of $\sigma (v)$. } \cr
%&2'.   \text{  If $\sigma  : \calb \rightarrow \calp$ does split,  }   \cr
%&\text{   \quad     then let  $\sigma (u)$ be the $0$th child of $\sigma (v)$ with
%probability ${1\over 2}$,  } \cr
%&\text{   \quad     and let  $\sigma (u)$ be the $1$st child of $\sigma (v)$ with
%probability ${1\over 2}$.   }    }$$
%
%

Recall that there will be one parallelogram for each vertex $t\in \calb _{h (\calp ) }$.
Since $x>1$, for each $t$, there is at most one possible value for
$\sigma (t)$ in $\calp_{h (\calp ) }$ such that $(x,y) \in P_{t, \sigma}$.  If such a
slope exists, call it
$S_{(x,y)} (t)$; if it does not exist, we will say $S_{(x,y)} (t) = \infty$.
%In other
%words, if $(x,y) \in P_{t, \sigma}$ for some choice $\sigma$, we require that if $u$ is
%an ancestor of $t$ then $\sigma(u)$ is an ancestor of $S_{(x,y)} (t)$.  In particular,
%this means that for every splitting vertex $v_1, ..., v_N$ that is an ancestor of
%$S_{(x,y)} (t)$, if $u$ is an ancestor of $t$ such that $\sigma (u) = v_j$, and if $w$ is
%the child of $u$ that is an ancestor of $t$, then $\sigma (w)$ must be the child of $v_j$
%that is the ancestor of $S_{(x,y)} (t)$.
The set of $t\in \calb _ {h (\calp) } $
for which $S_{(x,y)} (t)<\infty,$ call it the \it possible set \rm of $(x,y)$,
Poss$(x,y)$, will have at most $2^{N }$ elements.  Given a set of vertices
$V\subseteq \calb$, we denote by $<V>$ the tree generated by $V$, i.e., the subtree of
$\calb$ consisting of $V$ and all the ancestors of elements in $V$  (and all the edges
connecting these vertices).  Now consider the subtree $<Poss(x,y)> \subseteq \calb$.
Given $t\in Poss(x,y)$, there are at least $N$ ancestors of $t$, say $t_1, ..., t_N \supseteq t$ such that $t_j$
is an ancestor of $t_{j+1}$ and such that $\sigma (t_j) $ is a splitting vertex in
$\calp$.
Call such vertices \it choosing vertices \rm.  

Now let $\calc$ be the tree formed by all
the choosing vertices of $<Poss(x,y)>$, and edges connecting any pair of choosing vertices
$u,v \in <Poss (x,y)>$ such that $u \subseteq v$ with no choosing vertex $w$ such
that $u \subsetneqq w \subsetneqq v$.  Similarly, let $\calb ^* _N$ be the tree formed by all the splitting vertices of $\calp$,
with edges
connecting all the splitting vertices $u,v \in \calp$ such that $u \subseteq v$ and
there is no splitting vertex $w$ such that $u \subsetneqq w \subsetneqq v$.  Now $\calb
^* _N$ is
the binary tree of height $N$, i.e.,
$$\calb ^* _N = \calb \cap \left( \bigcup _{k=1} ^N \calb _k \right) ,$$
and $\calc$ is a subtree of $\calb ^* _N$.

So now we construct the sticky maps $\sigma : \calb \rightarrow \calp $ randomly as
follows:  to each edge $e$ in $\calc$, assign a random
variable $r=r(e)$ that takes on the values $0$ and $1$ with probabilities ${1\over 2}$.
We will write $e_{v,u}$ to denote the edge connecting $v$ to one of its children $u$.
If $r(e_{v,u}) = l$, we set
$$\sigma  (u)  = c_l ( \sigma  (v) ), $$
where, again, we use $c_l (w)$ to denote the $l$th child of a vertex $w$.

Let $k\in \{ 0,1,..., N \} $.  Given $v \in \calb _k$,
for $j=0,1,..., k$, define $A_j (v)$ to be the ancestor of $v$ at height $j$.
So if $(x,y) \in \calp _{t,\sigma}$, and if
$$A_{j+1} (S_{(x,y)} (t) ) = c_{b_j (t) }  (A_{j} (S_{(x,y)} (t) )  ),$$
for some sequence $b_j$ of zeros and ones depending on $t$,
then we must have
$$ r(e_{A_j(t) , A_j (t+1) }) = b_j (t).$$

%So if $(x,y) \in P_{t, \sigma}$, and if $\bar{v_1} (t) \supseteq \bar{v_2}(t)
%\supseteq ...\supseteq \bar{v_N}(t)$ are the splitting ancestors of $S_{(x,y)} (t)$,
%and $v_{j+1}$ is the $k_j$th child of $v_j$, and if  $\bar{u_1} (t) \supseteq
%\bar{u_2}(t)  \supseteq ...\supseteq \bar{u_N}(t)$ are the choosing ancestors of $t$,
%then we require $r(e_{u_j,u_{j+1} } ) = k_j (t)$.

Similarly, if we are to have $(x,y)
\in K_{ \sigma}$, then we must find a $t\in Poss(x,y)$ such that
$r(e_{A_j(t) , A_j (t+1) }) = b_j (t)$ for all $j= 1,...,N$.  Since $r$ takes on each
value $0$
or $1$ with probability ${1 \over 2}$, this requirement is equivalent to the following:
if we remove each edge of $\calc$ with probability ${1\over 2}$, we require that a path
remains from the root to $\calc _N$.  In the probability literature, the probability of
this outcome is called the survival probability of
Bernoulli$({1\over 2} )$ percolation on $\calc$, which we discuss below.

% \head \S 4 Percolation on Trees \endhead

Given a tree $\calt \subseteq \calb$ of height $N$, remove each edge with probability
${1\over 2}$.   Denote by $P(\calt )$ the probability that a path remains from the origin to
$\calt _N$.  A convenient way to compute this quantity is to view $\calt$ as an
electrical circuit.  Accordingly, we define the resistance of the tree $\calt$ as
follows:  place the positive node of a battery at the root of $\calt$, then identify all
vertices in $\calt _N$, and place the negative node of the battery at this new vertex.
For each edge at distance $k$ from the root, place a resistor of strength $2^k$.  The
resistance of the tree $\calt$, call it $R(\calt)$, is defined to be the resistance of
this circuit.  The following result of R. Lyons relates the resistance of $\calt$ to the
survival probability of Bernoulli $({1\over 2} )$ percolation on $\calt$.  We state and
prove a special case to keep the paper self-contained.  For a more general result, see
[L].  For more about probability on trees, see [LP].  The proof given here is from [BK] and 
actually holds when $\calt$ is a subset of the ternary tree, but it works for our purposes since 
the binary tree is a subtree of the ternary tree.

\proclaim{Theorem 4.1 (Lyons)} We have that
$$P(\calt) \lesssim {1 \over 2 + R(\calt) }.$$
\endproclaim
\demo{Proof of Claim 3.1 B}
Assuming Theorem 4.1, it remains to show that the resistance of the tree $\calc$ is
$\gtrsim {1\over N}$.
First recall that $\calc$ is a subtree of the truncated binary tree
$\calb ^*  _N$, so $R(\calc) \geq R(\calb ^* _N )$.  Now to compute a lower bound for
$R(\calb ^* _N )$, connect all vertices at height $k$ by an ideal conductor to make one node
$V_k$.  (This only decreases the resistance of the circuit.)  Now there are $2^k$
connections between
$V_k ^*$ and $V_{k+1} ^*$, each with resistance $2^k$.  If $R_k$ is the resistance
between $V_k ^*$ and $V_{k+1} ^*$, then
$${1\over {R_k} } = \sum _{1} ^ {2^k} {1\over {2^k} } = 1, $$
so $R_k =1$ for all $k = 1,2,..., N$.  Summing over $k$ results in $R(\calb ^* _N ) \geq N$.
\qed
\enddemo

\demo{Proof of Theorem 4.1} We prove this by induction on $n$. Clearly it is true for
constant 2, when
$n=0$.
We assume up to $n-1$, we have
$$P(\calt) \leq {12 \over 2 + R(\calt ) }.$$

We observe that we may view
$\calt$ as the root together with up to 3 edges connected to 3 trees $\calt_1,\calt_2,$
and
$\calt_3$.
(If some of these trees are empty, we assign them probabilty zero and infinite
resistance.)
We denote
$$P(\calt)=P_j,$$
and
$$R(\calt)=R_j.$$
Then
we have the recursive formulae
$$P(\calt )={1 \over 2} (P_1 + P_2 + P_3) - {1 \over 4} (P_1P_2 + P_1P_3 + P_2P_3)
+{1 \over 8} P_1 P_2 P_3 \tag 4.1$$
and
$${1 \over R(\calt)}={1 \over 2+2R_1} + {1 \over 2 + 2R_2}
+{1 \over 2 + 2 R_3}. \tag 4.2$$

Now we break into two cases. In the first case, we have ${12 \over 2 + R_j} > 2$ for
some $j$. Then we have $R_j < 4$. This implies $R(\calt)<10$ which implies
${12 \over 2 + R(\calt)} > 1$, so that we certainly have
$$P(\calt) \leq {12 \over 2 + R(\calt)}.$$

We define
$$Q_j={12 \over 2 + R_j}.$$
We may assume each $Q_j \leq 2$.
Observe that if we define
$$F(x,y,z)= 1 -(1-{1 \over 2} x)(1-{1 \over 2}y)(1 -{1 \over 2}z),$$
on the domain $[0,2] \times [0,2] \times [0,2]$ then $F$ is monotone increasing in
each variable.
Therefore we have that
$$\eqalign{ P(\calt) &=  F(P_1,P_2,P_3) \cr
                         &\leq  F(Q_1,Q_2,Q_3) \cr
                         &\leq {1 \over 2}(Q_1+Q_2+Q_3) - {1 \over 6} (Q_1Q_2 + Q_1
Q_3 +Q_2 Q_3) }.\tag 4.3$$
Note that the equality is (4.1), while for the two inequalities we have used that the
$Q$'s are
$\leq 2$.

Now plugging into (4.3), the definition of the $Q$'s,  we obtain
$$\eqalign{ P(\calt )
&\leq {12 \over 2} [ {(R_1+2)(R_2+2) + (R_1+2)(R_3+2) + (R_2+2)(R_3+2) - {12  \over 6}
(R_1+R_2+R_3 + 6) \over (R_1+2)(R_2+2)(R_3+2) }] \cr
&\leq {12\over 2} [ {(R_1+2)(R_2+2) + (R_1+2)(R_3+2) + (R_2+2)(R_3+2) - {12 \over 6}
(R_1+R_2+R_3 + 6) \over (R_1+2)(R_2+2)(R_3+2) - 4R_1-4R_2-4R_3-13}] \cr
&\leq {12 \over 2} [ {(R_1+1)(R_2+1) + (R_1+1)(R_3+1) + (R_2+1)(R_3+1)
\over (R_1+2)(R_2+2)(R_3+2) -4R_1-4R_2-4R_3-13}] \cr
&= {12 \over R(\calt) +2}. } $$
Here the second inequality is by
decreasing the denominator and the third inequality is by increasing the numerator.

\qed \enddemo

 \head \S 5  The Lacunary Case \endhead

To complete the proof of Theorem 0.1, it remains to show the following proposition:
\proclaim{Proposition 5.1}  If $\calt _{\Omega} $ is lacunary of order $N$, then there exists a constant $C$ 
such that
$$|| M _ { \Omega } f ||_p \leq CN ||f||_p .$$
\endproclaim 
\proclaim{Remark 5.2}
As noted earlier, if $\tom$ has splitting number $N$, then $\tom$ is lacunary of order $N$, and hence 
$$ || M _{\Omega} f ||_p \leq C N ||f||_p . $$
\endproclaim  
As mentioned earlier, the result in Proposition 5.1 was published in [SS].  The proof given here extends ideas 
present in [NSW], and follows Alfonseca [A].

Recall that each ray $R\in \tom$ is identified with a real number $\alpha (R) \in [0,1]$.  
If $\calt$ is a tree, we will define
$$ \alpha ( \calt ) = \{ \alpha (R) \colon R \in \partial T\} .$$

Now write $v_{\theta}$ to denote the unit vector with slope $\theta$, and define 
$$ M_{\theta} f(x) = \sup _{h>0} {1 \over {2h} } \int _{-h} ^h f(x+v_{\theta} t) dt .$$  
%If $\theta = \alpha (R)$, then define $v_{\theta} = v_R$ and $M_{\theta} f= M_R f$.
If a tree $\mathcal L $ is lacunary of order $1$, then there is a ray, called $\lim \mathcal L $ as in Remark 1.1, such that 
every splitting vertex in $\mathcal L $ lies along $\lim \mathcal L $.  
Also define $\beta _j (\mathcal L )$ to be the (unique, if it exists) element of $\alpha (\mathcal L) $ such that 
$d(\alpha (\lim \Cal L ) , \beta _j ) = 2^{-j}$, where again $d$ denotes the dyadic distance on real numbers, and let 
$$\Omega _j (\mathcal L ) := \{ \beta \in  \alpha (\mathcal  L ) \colon d( \alpha (\lim \mathcal L ) , \beta ) = 2^{-j} \}.$$
\proclaim{Proposition 5.3}
Fix $1<p\leq \infty$. Let $\Omega ^* \subseteq \Omega $.  If there exists a lacunary tree $\mathcal L$ of order $1$ such that 
$\mathcal L = \calt _{\Omega ^* }$, then there exists a constant $C$, depending only on $p$, such that
$$||M_{\Omega} f||_p \leq C ||f||_p \left( 1 + \sup _j ||M _{\Omega _j} || _{L^p \rightarrow L^p} \right).$$
\endproclaim
A simple iteration of Proposition 5.3 will give us Proposition 5.1:  For if $\tom$ is lacunary of order $N$, then 
there exists $\Omega ^* \subseteq \Omega$ and a lacunary tree $\mathcal L$ of order $1$ such that 
$\calt _{\Omega ^* } = \mathcal L$.  Since $\tom$ is lacunary of order $N$, we may actually choose such an $\mathcal L$ 
so that $\calt  _{\Omega _j (\mathcal L ) } $ is lacunary of order $N-1$ for all $j$.  But then we may apply the propostion 
again to the sets $\Omega _j$ and repeat to get $||M_{\Omega} f || _p \leq CN || f || _p $.  It remains to prove Proposition 5.3.

\demo{Proof of Proposition 5.3}
For convenience, we will rotate the plane so that $\alpha (\lim \mathcal L ) =0$.  
Let $\delta _j = { {12}\over {20} } 2^{-j}$.  Let $A_j$ be an interval centered around the dyadic interval containing 
the sets $\Omega _j$ such that $|A_j| = \delta _j $ and $dist(A_j ^c , \Omega _j ) \geq {1\over {50} } 2^{-j} $.  Also let 
$\widetilde{A_j} = {{13}\over {11}} A_j$ and $\widetilde{ \widetilde{A_j} } = {{14}\over {11}} A_j$.  Then define 
$$ \eqalign{   \Delta _j &= \{ (x,y)\in \Bbb R ^2 \colon {y\over x} \in A_j \} \cr
	\widetilde{ \Delta _j } &= \{ (x,y)\in \Bbb R ^2 \colon {y\over x} \in \widetilde{A_j } \} \cr
	\widetilde{ \widetilde{ \Delta _j } } &= \{ (x,y)\in \Bbb R ^2 \colon {y\over x} \in \widetilde{ \widetilde{ A_j } } \} ,
}$$
and let $\omega_j$ be a $C^{\infty}$ function away from the origin, homogeneous of degree zero, such that 
$\omega_j \equiv 1 $ on $ \Delta _j$ and $\omega_j \equiv 0 $ outside $ \widetilde{ \Delta _j } $.  Similarly, define 
$\widetilde{\omega_j } $ with respect to the sectors $\widetilde{ \Delta _j }$ and $\widetilde{ \widetilde{ \Delta _j } }$.  Define 
$$ \widehat{ S_j f } = \omega _j \widehat{f}, \quad \quad \quad  \widehat{ \widetilde{S_j f } } = \widetilde{\omega _j } \widehat{f} .$$

Now let $\psi : \Bbb R \rightarrow \Bbb R$, $\psi \in C^{\infty}$, $\psi \geq 0$, be such that 
$\psi \equiv 1$ on $[-1,1]$ and $\psi \equiv 0$ outside $[-2,2]$.  For $j\geq 1$ and $\theta \in \Omega _j$, define 
$$ N_{h,j,\theta} f (x) = {1\over {2h}} \int_{-\infty} ^{\infty} \psi \left( {t\over h} \right) f(x + v_{\theta} t ) dt .$$
We will only consider nonnegative $f$, and for such $f$, $\sup _{h>0} N_{h,j,\theta} f (x) \sim M_{\theta} f(x)$.  
Let $m= \widehat \psi $ and let $\phi : \Bbb R ^2 \rightarrow \Bbb R$ be $ C^{\infty} _c $ with $\phi (\xi ) =1$ when $|\xi | \leq 1$.  
Now we write 
$$ \eqalign{   \widehat{ N_{h,j,\theta} f } (\xi ) & = m(h\xi _1 + h \xi _2 \theta ) \widehat{f} (\xi ) \cr
& = \phi(h\delta _j \xi )  m(h\xi _1 + h \xi _2 \theta ) \widehat{f} (\xi ) \cr 
& + (1- \phi(h\delta _j \xi ) ) ( 1- \omega _j (h\xi ) )    m(h\xi _1 + h \xi _2 \theta ) \widehat{f} (\xi ) \cr
& + (1- \phi(h\delta _j \xi ) ) \omega _j (h\xi )  m(h\xi _1 + h \xi _2 \theta ) \widehat{f} (\xi )   \cr
& =: I _{h,j,\theta} + II _{h,j,\theta} + III _{h,j,\theta} .    }  $$
To control $I _{h,j,\theta}$, we note that $\phi$ is Schwartz, and hence $I _{h,j,\theta}$ is controlled by the strong maximal 
function $M ^ {S} _ {\beta _j} $ with respect to the axes with slopes $\beta_j$ and $\beta _j + \pi $, with constants independent 
of $j$.  The term $II _{h,j,\theta}$ can be estimated in the same way, giving us 
$$ \sup _h \sup _{\theta \in \Omega _j } I _{h,j,\theta} f(x) + II _{h,j,\theta} f (x) 
\leq C M _{\beta _j } ^S f (x) .   \tag { $*$  }  $$
Hence 
$$ || \sup _{j,h,\theta \in \Omega _j}  (I _{h,j,\theta} + II _{h,j,\theta}  ) f ||_p 
\leq C || M_{\Omega ^* } f ||_p \leq C || f || _p ,    \tag { $**$  }   $$ 
by the result in [NSW], since $T_{\Omega ^* }$ is lacunary of order $1$.  
It remains to control $III _{h,j,\theta}$.

We will assume, for now, that $\Omega$ is a finite set, and obtain a bound independent of the size of $\Omega$.  
Since $\Omega $ is finite, there is some minimal constant $C(\Omega)$ such that for $f\in L^p$,
$$ \left|\left|  \sup _{j,h,\theta \in \Omega _j}  |  III _{h,j,\theta}  f | \right| \right|  _p \leq C(\Omega) ||f||_p .$$
Note that $N_{h,j,\theta} (g) \leq N_{h,j,\theta} (|g|) $ for any $g,h,j, \theta$, and let $\{ g_j \} $ be a sequence of functions.  
Then by $(*)$, $(**)$ and the decomposition of $N_{h,j,\theta}$, we have 
$$ \eqalign{  \left| \left|  \sup _{j,h,\theta \in \Omega _j}  |  III _{h,j,\theta}  (g_j) | \right| \right|_p  
& \leq   \left|\left|  \sup _{j,h,\theta \in \Omega _j}  |  III _{h,j,\theta} ( \sup _{j} | g_j |    ) | \right| \right|_p  + c_0 \left| \left|  \sup _{j} | g_j |  \right| \right|_p  \cr
& \leq (C(\Omega ) + c_0 ) \left| \left|  \sup _{j} | g_j |  \right| \right| _p  . } $$
In fact, if $C(\Omega) \leq c_0$  independent of $\Omega$, we are already finished with Proposition 5.3, 
so we can assume otherwise and estimate the quantity above by $2C(\Omega ) ||f||_p $.  In addition , it is clear that 
$$ \left|\left| \left(  \sum _j \sup _{h,\theta \in \Omega _j}  |  III _{h,j,\theta}  (g_j) | ^p \right) ^ {  {1\over p} } \right| \right| _p 
\leq \left ( \sup _j || \sup _{h,\theta \in \Omega _j}    III _{h,j,\theta}  || _ {L^p \rightarrow L^p } \right) 
\left|\left| \left( \sum _j | g_j | ^p \right ) ^{  {1\over p} }\right| \right| _p .$$
Interpolating yields 
$$ \left|\left| \left(  \sum _j \sup _{h,\theta \in \Omega _j}  |  III _{h,j,\theta}  (g_j) | ^2 \right) ^ {  {1\over 2} } \right| \right| _p 
\lesssim C(\Omega ) ^{1-{p \over 2} } 
\left ( \sup _j || \sup _{h,\theta \in \Omega _j}    III _{h,j,\theta}  || _ {L^p \rightarrow L^p } \right) ^{p\over 2}  
\left|\left| \left( \sum _j | g_j | ^2 \right ) ^{  {1\over 2} }\right| \right| _p .$$
Recall that $ III _{h,j,\theta} f $ has frequency support in $\widetilde{ \Delta _j }$ so 
$$ \eqalign{   \left| \left | \sup _{j,h,\theta \in \Omega _j}  | III _{h,j,\theta} f | \right| \right|  _p 
&  \leq \left| \left| \left( \sum _j \sup _{h,\theta \in \Omega _j}  
| III _{h,j,\theta} ( \widetilde{ S _j } f ) | ^2 \right) ^{1\over 2} \right| \right| _p    \cr
& \lesssim C(\Omega ) ^{1-{p\over 2} }  
\left( \sup _j || \sup _{h,\theta \in \Omega _j}    III _{h,j,\theta}  || _ {L^p \rightarrow L^p } \right) ^{p\over 2}  
\left| \left| \left( \sum _j | \widetilde{ S _j } f | ^2 \right ) ^{  {1\over 2} }\right| \right| _p .  }  $$
But $C(\Omega)$ is minimal, and one can show
$\left|\left| \left( \sum _j | \widetilde{ S _j } f | ^2 \right ) ^{  {1\over 2} } \right| \right| _p \leq C ||f||_p $
 by using Rademacher functions and the Marcinkiewicz multiplier theorem, as in [NSW], so
$$ C(\Omega) \lesssim C(\Omega) ^{1-{p\over 2} }  
\left( \sup _j || \sup _{h,\theta \in \Omega _j}    III _{h,j,\theta}  || _ {L^p \rightarrow L^p } \right) ^{p\over 2} ,$$
and hence
$$ C(\Omega)  \lesssim   \sup _j || \sup _{j,h,\theta \in \Omega _j}    III _{h,j,\theta}  || _ {L^p \rightarrow L^p }  . $$
However, 
$$ \eqalign{ III _{h,j,\theta} (f) & \lesssim N_{h,j,\theta} ( | \check{\tilde{\omega _j } } \ast f | )  \cr
& \lesssim M_{\Omega _j} (  | \check{\tilde{\omega _j } } \ast f | ) \cr
& \lesssim M_{\Omega _j} \left( ( \sum _j | \widetilde{ S _j } f | ^2 ) ^{  {1\over 2} } \right ) , }$$
so 
$$  \sup _j \left| \left| \sup _{h,\theta \in \Omega _j}    III _{h,j,\theta}  \right| \right| _ {L^p \rightarrow L^p }  
\lesssim \sup _j \left| \left| M_{\Omega _j} \right| \right| _{L^p \rightarrow L^p} ,$$
and this proves Proposition 5.2. \qed
\enddemo

\end